\documentclass[12pt]{article}
\usepackage[mathscr]{eucal}
\usepackage{amsbsy}
\usepackage{amssymb}
\usepackage{amsmath}
\usepackage{amsthm}
\usepackage{graphics}
\usepackage{epsfig}

\usepackage{latexsym}

\newcommand{\dd}{\,\mathrm{d}} % differential element in integrals
\newcommand{\E}{\mathbf{E}\,}

\def\phi{\varphi}

\def\Bl1{{\bf 1}}
\def\B2{{\bf 2}}
\def\B0{{\bf 0}}

\def\a{\alpha}
\def\b{\beta}
\def\=A8{\"o}

\newcommand{\beq}{\begin{equation}}
\newcommand{\eeq}{\end{equation}}
\newcommand\beqn{\begin{displaymath}}  % no number
\newcommand\eeqn{\end{displaymath}}

\theoremstyle{plain}
\newtheorem{teo}{Theorem}
\newtheorem{lem}[teo]{Lemma}

\theoremstyle{definition}

\newtheorem{remark}[teo]{Remark}

\begin{document}

\title{Once more on comparison of tail index estimators
\footnotemark[0]\footnotetext[0]{ \textit{Short title:}  Comparison of tail index estimators}
\footnotemark[0]\footnotetext[0]{%
\textit{MSC 2000 subject classifications}. Primary 62F12, secondary
62G32, 60F05 .} \footnotemark[0]
\footnotemark[0]\footnotetext[0]{ \textit{Corresponding
author:} Vygantas Paulauskas, Department of Mathematics and
 Informatics, Vilnius university, Naugarduko 24, Vilnius 03225, Lithuania,
 e-mail:vygantas.paulauskas@mif.vu.lt}}

%\footnote{The research was supported by the bilateral
%France-Lithuania scientific project Gilibert  and Lithuanian State
%Science and Studies foundation (V-07058).} }

\author{ Vygantas Paulauskas$^{\text{\small 1}, \text{\small 2}}$ and Marijus Vai{\v c}iulis$^{\text{\small 2}}$ \\
{\small $^{\text{1}}$ Vilnius University, Department of Mathematics
and
 Informatics,}\\
{\small $^{\text{2}}$
  VU Institute of Mathematics and Informatics,
 Vilnius}
 }

%\date{}

%\begin{document}

\maketitle

\begin{abstract}
We consider heavy-tailed distributions and compare the well-known
estimators of the tail index, based on extreme value theory with a
comparatively recent estimator based on a different idea.
\end{abstract}
\vfill
\eject
\section{Introduction}% and formulation of results}

During last several decades it was demonstrated that in many fields
of applied probability the so-called heavy-tailed distributions play
an important role. One of the main problems, connected with
heavy-tailed distributions is the estimation of the tail index - a
parameter, which characterizes the heaviness of the tail of a
distribution. The problem can be formulated as follows. Let us
consider a sample $X_1, \dots , X_N$ of size $N$ taken from a
heavy-tailed distribution function (d.f.) $F$, that is, we assume
that $X_1, \dots , X_N$ are independent identically distributed
(i.i.d.) random variables with a d.f. $F$ satisfying the following
relation for large $x$:
\begin{equation}\label{cond1}
1-F(x)=x^{-\a}L(x).
\end{equation}
Here $\a >0$, \  $L(x) >0$ for all $x>0$ and$L$ is a slowly varying
 at infinity function:
$$
\lim_{x\to \infty}{L(tx) \over L(x)}=1.
$$
If we have only condition (\ref{cond1}) without any additional
information about the function $L$, it is difficult to get good
properties, such as the asymptotic normality, of  an estimator of
the parameter $\a$. Therefore the main stream of papers dealing with
the tail index estimation uses the so-called second-order condition
of regular variation . During last years even the third order
condition on $F$ was introduced (see, for example \textsc{Fraga
Alves} et al(2006)). In our paper, like in \textsc{Paulauskas}
(2003), we shall use the second order condition in the form of the
relation (\ref{secorder})  below.

 Let $X_{N,1}\le X_{N,2}\le \cdots \le
X_{N,N}$ denote the ordered statistics of $X_1, \dots ,X_N.$ Most of
tail index estimators are based on ordered statistics and estimate
the parameter $\gamma=1/\alpha$. One of the most popular  estimators
to estimate the parameter $\gamma=1/\alpha$ was proposed by \textsc{Hill} (1975):
\[
\gamma_{N,k}^{(1)}={1\over k}\sum_{i=0}^{k-1}\log X_{N,N-i} - \log
X_{N,N-k},
 \]
where $k$ is some number satisfying $1\le k\le N.$ We also list some
other estimators based on ordered statistics:
\begin{eqnarray*}
&&\hspace{-1pc}\gamma_{N,k}^{(2)}=(\log 2)^{-1}\log {{X_{N,N-[k/4]}
- X_{N,N-[k/2]}}\over
{X_{N,N-[k/2]} - X_{N,N-k}}}, \\
&&\hspace{-1pc} \gamma_{N,k}^{(3)}=\gamma_{N,k}^{(1)}+1-\textstyle
\frac 12\big
(1-(\gamma_{N,k}^{(1)})^2/M_N\big )^{-1},\\
&&\hspace{-1pc}\gamma_{N,k}^{(4)}={M_N \over {2\gamma_{N,k}^{(1)}}},
\end{eqnarray*}
where
\[
M_N={1\over k}\sum_{i=0}^{k-1}(\log X_{N,N-i} - \log X_{N,N-k})^2.
\]
The estimator $\gamma_{N,k}^{(2)} $ was proposed by \textsc{Pickands} (1975)
, $\gamma_{N,k}^{(3)}$ -- in \textsc{Dekkers} et al (1989)
 and
$\gamma_{N,k}^{(4)}$ -- by C.~G. de Vries (see e.g. \textsc{de Haan} and  \textsc{Peng} (1998)). There are
many papers devoted to the modifications of  the estimators
$\gamma_{N,k}^{(i)}, \ i=1,2,3,4$. Let us mention several of them (in chronological order):
\textsc{Weissman} (1978), \textsc{Smith} (1987), \textsc{Resnick} and \textsc{Starica} (1997), \textsc{Geluk }and \textsc{Peng} (2000),
\textsc{Fraga Alves} (2001),
\textsc{Gomes} and \textsc{Martins} (2002), \textsc{Fraga Alves} et al (2003), \textsc{Li} et al (2008), \textsc{Gomes} et al (2008).

All  estimators $\gamma_{N,k}^{(i)}, \ i=1,2,3,4$ contain one
additional parameter $k$, which has clear intuitive meaning in the
case of all above written estimators $\gamma_{N,k}^{(i)}$: how many
the largest values from the ordered statistics must be taken in
order to have good properties (consistency, asymptotic normality,
etc) of the estimator.  We presented these four estimators because
in \textsc{de Haan} and \textsc{Peng} (1998) all these estimators
are compared and it is shown that none of these estimators dominates
the others. It turned out that for different values of the
parameters $\gamma$ and~$\rho$ (the parameters characterizing the
so-called second-order asymptotic behavior of~$F$, which will be
introduced below) different estimators have the smallest asymptotic
mean-squared error.

In  \textsc{Davydov} et al (2000) (see also \textsc{Davydov} and
\textsc{Paulauskas} (1999)) there was proposed a new estimator,
based on a different idea, which came when considering rather
abstract objects -- random compact convex stable sets. Although
originally in \textsc{Davydov} et al (2000) an estimator was
constructed to estimate the index of a multivariate stable
distribution ( even there was the restriction $0<\a <1$ for this
index) and the main tool in the proof was the relation between
exponential distribution and ordered statistics, in
\textsc{Paulauskas} (2003) it was noted that the same construction
of the estimator can be based on a different idea and that this idea
 can be employed in the context of a general tail index estimation.
 We shall return to this point after introducing this new estimator.
%Using the second-order condition for the tail of a distribution main
%properties of the estimator were established.
The construction of the estimator is as follows.

We divide the sample into $n$ groups $V_1, \dots, V_n$, each group
containing $m$ random variables, that is, we assume that $N=n\cdot
m$ and $V_i=\{X_{(i-1)m+1}, \dots , X_{im+1} \}$. (In practice at
first $m$ is chosen and then $n=[N/m]$ is taken, where $[x]$ stands
for the integer part  of a number $x>0.$) Let
\[
M_{ni}^{(1)}=\max \{ X_j{:}\ X_j\in V_i\}
\]
and let $M_{ni}^{(2)}$ denote the second largest element in the same
group $V_i.$ Let us denote
\begin{equation} \label{definitions}
\kappa_{ni}={M_{ni}^{(2)}\over M_{ni}^{(1)}}, \qquad
S_n=\sum_{i=0}^n\kappa_{ni},\qquad
 Z_n=n^{-1}S_n.
\end{equation}

From now instead of (\ref{cond1}) we require stronger condition. Let
us assume that we have a sample $X_1, \dots ,X_N$ from distribution
$F$, which satisfies the second-order asymptotic relation (as $x\to
\infty$) \beq \label{secorder}
1-F(x)=C_1x^{-\alpha}+C_2x^{-\beta}+\mathrm{o}(x^{-\beta }),
%\eqno(2)
\eeq with some parameters $0<\alpha<\beta\le \infty.$ It seems that
Hall (see \textsc{Hall} (1982)) was the first who considered this
condition in the context of tail index estimation and, under this
condition, proved asymptotic normality of Hill estimator. Note that
the case $\beta =\infty$ corresponds to Pareto distribution,
$\beta=2\alpha$ -- to stable distribution with exponent
$0<\alpha<2.$ In \textsc{de Haan} and  \textsc{Peng} (1998) (and in
many other papers dealing with tail index estimation as well) more
general second-order asymptotic relation is used in a different form
with parameters $\gamma=1/\alpha$ and $\rho=\alpha-\beta$ and in a
more general context of the extreme-value index, when the parameter
$\gamma$ can take negative values, too. Namely, let
 $U$ denotes the right continuous
inverse of the function $1/(1-F)$. Suppose that there exits a
function $A(t)$, which ultimately has constant sign and tends to
zero as $t \to \infty$, then the relation
 \beq \label{gensecor}
\lim_{x \to \infty} \frac{\frac{U(tx)}{U(t)}-x^{\gamma}}{A(t)} =
x^{\gamma} \frac{x^{\rho} -1}{\rho}, \quad \rho \le 0, \eeq
 serves for the
definition of the second order regularly varying tail $1-F$.

As it was mentioned above, the estimator $Z_n$ from
(\ref{definitions}) (in a different context of a sample from
multivariate stable distribution) was based on the following
relation (see  \textsc{LePage} et al (1981))
$$
\big (M_{ni}^{(1)}, M_{ni}^{(2)}\big )m^{-1/\alpha} \buildrel \rm D
\over \longrightarrow_{N\to \infty} (\Gamma_1^{-1/\alpha},
\Gamma_2^{-1/\alpha}),
$$
where $\Gamma_i=\sum_{j=1}^i \lambda_j$ and $\lambda_j, j\ge 1,$ are
i.i.d. standard exponential random variables  and $\buildrel \rm D
\over \longrightarrow$ denotes convergence in distribution, and the
fact that
$$
\E\Big (\frac{\lambda_1}{\lambda_1+\lambda_2}\Big
)^{1/\alpha}=\frac{\a}{1+\a}.
$$
 Here may be it is worth to mention that most estimation of tail
index procedures, starting from \textsc{Hill}(1975) and
\textsc{Pickands} (1975) papers, are based on the relation between
order statistics and exponential distributions (the Renyi
representation theorem): if $X_1, \dots , X_n$ is a sample from a
continuous strictly increasing d.f. $F, \ F(0)=0,$ and $X^{(1)}\ge
\dots \ge X^{(n)}$ is the order statistics, then
$$
X^{(i)}=F^{-1}\left(\exp\left\{-\sum_{j=1}^i\lambda_j (n-j+1)^{-1} \right\}\right), \quad i=1,2, \dots ,n
$$
and
$$
\lambda_j=(n-j+1)\big (\ln F(X^{(j-1)})- \ln F(X^{(j)})\big ),\quad
j=1,2, \dots ,n
$$
In \textsc{Paulauskas} (2003) it was noted that estimator from
(\ref{definitions}) can be based on a different idea. If we take two
independent random variables $X$ and $Y$ with the same Pareto
distribution
$$
F(x)=1-C_1x^{-\alpha}, \quad x\ge C_1^{1/\a},
$$
and denote
\begin{equation} \label{Pareto}
W={{\min (X, Y)}\over {\max (X, Y)}},
\end{equation}
then it is not difficult to verify that, denoting $p=\a/(1+\a),$ we
have $\E W=p$ (since $W$ is invariant under scale transformation, we
can take $C_1=1$ and in the sequel we shall refer to that case as a
standard Pareto distribution).
 Therefore in the case of the Pareto distribution
 quantity $Z_n$, as an estimator for the parameter
$p$ (we shall denote it by $\hat p$ and as in \textsc{Qi} (2010) we
shall call it as DPR estimator), is nothing but the sample mean for
a bounded random variable, moreover, in this case the best choice is
to take m=2. If the underlying distribution $F$ is not Pareto, but
satisfies (\ref{secorder}), then it is natural to expect that for
large $m$ \ $\E\hat p=\E\kappa_{n1}$ will be close to $p.$ In
\textsc{Paulauskas} (2003) the following estimate , which was the
main ingredient  in the proof of the asymptotic normality of the
estimator $\hat p$, was given (see Lemma in \textsc{Paulauskas}
(2003)) \beq \label{gamma_m}|\gamma_m| \le C_0 m^{-\zeta}, \eeq
where $\gamma_m=\E \hat{p}-p,$ \ $\zeta=(\beta-\a)/ \a$, \ and $C_0$
is a constant depending on $C_1, C_2, \a,$ and $\beta$. The DPR
estimator  is a sum of i.i.d. random variables, therefore $\E\hat
p=\E\kappa_{n1}$ and $\gamma_m$ gives the bias of our estimator.
Based on  (\ref{gamma_m}) the following result was proved in
\textsc{Paulauskas} (2003).
\begin{teo}\label{thm1}
Let us suppose that $F$ satisfies $(2)$ with $0<\a<\beta< \infty. $
If we choose
\[
n=\varepsilon_N N^{2\zeta/(1+2\zeta)},\qquad
  m=\varepsilon_N^{-1} N^{1/(1+2\zeta)},
\]
where $\varepsilon_N \to 0,$ as $N\to \infty$ , then
\beq \label{clt_paul}
\sqrt n
(\hat p- p)\buildrel \rm D \over \longrightarrow_{N\to \infty}
N(0,\sigma^2),
%\eqno(8)
\eeq where $\sigma^2=\lim_{n\to \infty}\sigma_n^2={\alpha
((\alpha+1)^2(\alpha+2))^{-1}}$.
\end{teo}

 Now we can give the exact asymptotic behavior of
$\gamma_m$, this allows us to choose $m$ in an optimal way and to
compare DPR estimator  with these four estimators, listed above, in
the same manner as it was done in \textsc{de Haan} and \textsc{Peng}
(1998) (therefore in the title of the paper there are words ``once
more"). Our main result can be formulated as follows. (We write $a_n
\sim b_n$ if $\lim_{n\to \infty} a_nb_n^{-1}=1$. )
\begin{teo}\label{thm2}
Let us suppose that $F$ satisfies $(2)$ with $0<\a<\beta < \infty $
and $C_1>0$. Then we have
\begin{equation} \label{bias_ratio_asymp}
\gamma_m \sim \chi m^{-\zeta}, \quad (m \to \infty),
\end{equation}
 where
$$
\chi=\chi(C_1, C_2,\alpha, \beta)=\frac{C_2 \beta \zeta
\Gamma(\zeta+1)}{C_1^{\zeta+1} (\alpha+1)(\beta+1)}.
$$
For sufficiently large $N$ (ensuring that $m_{opt}\ge 2$) taking
\begin{equation}\label{opt4}
    m_{opt}=\left( \frac{2 \zeta \chi^2 }{\sigma^2}\right)^{1/(1+2\zeta)} N^{1/(1+2\zeta)}
\end{equation}
we get that MSE (mean square error) is minimal
\begin{equation}\label{opt5}
 \E\left(\hat{p} -p \right)^2 \sim \left(1+2\zeta\right) \left(
\frac{\chi^2  \sigma^{4 \zeta}}{ (2\zeta)^{2\zeta}
N^{2\zeta}}\right)^{1/(1+2\zeta)}.
\end{equation}
Under this choice of $m$ we have asymptotic normality
 \beq \label{clt_eq}
\sqrt n (\hat p- p)\buildrel \rm D \over \longrightarrow_{N\to
\infty} N(\mu,\sigma^2), \eeq where $\sigma^2$ is the same as in
Theorem \ref{thm1}  and $\mu=\sigma (2 \zeta)^{-1/2}{\rm  \ sgn}(\chi)$.
\end{teo}

\begin{remark} The estimator $\hat p$ (as all other introduced above tail index
estimators) is invariant with respect to scale transformation, while
condition (\ref{secorder}) is not: if a random variable $X_1$
satisfy this condition, a random variable $AX_1$, where $A>0$ with
distribution function $F_A$ satisfies the relation
$$
1-F_A(x)=C_1A^\alpha x^{-\alpha}+C_2A^\beta
x^{-\beta}+\mathrm{o}(x^{-\beta }),
$$
therefore the constants in the second order relation are not
invariant. But it is easy to see that the ratio $C_1^\beta/
C_2^\alpha $ is invariant and all quantities in relations
(\ref{bias_ratio_asymp}), (\ref{opt4}), and (\ref{opt5}) depend
exactly only on this ratio.
\end{remark}

As it was mentioned, this result allows to compare the estimator
$\hat p$ with four estimators mentioned above, and this is done in
Section 3. Although according to the chosen criteria Hill estimator
$\gamma_{N,k}^{(1)}$ and estimator $\gamma_{N,k}^{(4)}$
asymptotically perform better than the estimator $\hat{p}$, relation
between  other two estimators and $\hat{p}$ is the same as in paper
\textsc{de Haan} and \textsc{Peng} (1998): for some values of
parameters $\a, \beta$ estimator $\hat{p}$ performs better than
$\gamma_{N,k}^{(2)}$ and $\gamma_{N,k}^{(3)}$. But here it is worth
to mention  the simple structure of the DPR estimator  and the fact
that it is well suited for recursive calculations (for example, when
we have the so-called tick-by-tick financial data and we need tail
index estimation in real time), see the monograph \textsc{Markovich}
(2007). There are situations (such examples are mentioned in
\textsc{Oi}(2010)) when data can be divided naturally into blocks
but only few of largest observations within  blocks are available.
In such situations estimator $\hat{p}$ can be applied while all
other mentioned estimators are not applicable. Also the estimator
$\hat{p}$ is well adapted for detecting a change in tail index, see
a paper \textsc{Gadeikis} and \textsc{Paulauskas} (2005), where
estimator  $\hat{p}$ was used to analyze financial crisis in Asian
markets in 1997-98 and the results were compared with analogous
analysis using Hill estimator in \textsc{Quintos} et al(2001). And
the main factor why we think that DPR estimator deserves the
attention of statisticians  is the possibility of several promising
modifications of the estimator. One such modification is to
introduce an additional parameter $r>0$ and to consider the
estimator
$$
{\hat p}_r=n^{-1}S_{n,r}, \quad S_{n,r}=\sum_{i=1}^n\kappa_{ni}^r.
$$
Again, using standard Pareto distribution, it is easy to calculate
that ${\hat p}_r$ estimates the quantity $\a(r+\a)^{-1}$. As a
matter of fact, when preparing  the paper \textsc{Paulauskas} (2003)
the first named author had considered the estimator ${\hat p}_2$
(which to some extent resembles the quantity $M_N$, see definition
of estimators $\gamma_{N,k}^{(3)}$ and $\gamma_{N,k}^{(4)}$ ), but
realized that there is no gain in changing the first moment by the
second one. Now it turns out that it is worth to take $0<r<1,$ and
we are able to prove asymptotic normality (under appropriate
assumptions on $m$) for a fixed $r$. Also we can show that between
two estimators with fixed parameters $0<r'<r''\le 1$, the smaller
asymptotic MSE has ${\hat p}_{r'}$. Unfortunately, at present we do
not know how to choose optimally $r$, which in general can be
dependent on $\a, \b,$ and even $N$.

Let us consider general construction of  modifications of the DPR
estimator.  Take some function $f: [0,1]\to [0,\infty]$ such that
$\E f(W)$ exists where $W$ is from (\ref{Pareto}), then this
expectation will be some function  of $\a$ and, of course, on
function $f$. Let us denote this function by $h_f(\a)$, that is
$h_f(\a)=\E f(W)$. If $h_f$ is a one-to-one map from $[a,b]$ to
$[c,d]$ with $[a,b]$ and $[c,d]$ being subsets of $[0,\infty]$, then
estimating the quantity $h_f(\a)$ and taking the inverse function we
get an estimator for $\a$ (with the restriction $a\le \a \le b$ if
$0<a<b<\infty$). Therefore it is natural to consider statistic of
the form
\begin{equation}\label{gDPR}
 \frac{1}{n}  \sum_{i=1}^nf(\kappa_{ni}),
\end{equation}
obtaining large class of modifications of the estimator ${\hat p}$,
developed in \textsc{Paulauskas} (2003). The estimator ${\hat p}$ is
obtained taking $f_1(x)=x$, then $h_{f_1}(\a)=\a/(1+\a).$ The above
mentioned modification  ${\hat p}_r$ is obtained by taking
$f_r(x)=x^r, \ r>0$ and $h_{f_r}(\a)=h_r(\a)=\a/(r+\a).$ Estimators
of the type (\ref{gDPR}) we shall call generalized DPR estimators,
in short GDPR.

In a recent paper \textsc{Qi} (2010)  one more estimator is
introduced, which can be considered as connecting ideas of DPR and
Hill estimators  At first the procedure is the same as for DPR
estimator - division of the sample in $n$ groups with $m$ elements
in each group. But then instead of taking two largest elements in
each group Qi takes  Hill estimator in each group, namely, taking
$s+1$ ($1\le s\le m-1$) largest values in each group, then averaging
them over groups and obtaining the following estimator of the
parameter $\gamma=\a^{-1}$:
\begin{equation}\label{Qi-s}
 \gamma_N(s)=\frac{1}{ns}  \sum_{i=1}^n \sum_{j=1}^s\big ( \log M_{ni}^{(j)} - \log
M_{ni}^{(s+1)}\big ),
\end{equation}
where $M_{ni}^{(1)}\ge \dots M_{ni}^{(m)}$ is ordered statistic from
$V_i.$ With $s=1$ estimator (\ref{Qi-s}) becomes of the form
(\ref{gDPR}) with $f_{\ell}(x)=-\log x.$ It is not difficult to
calculate that for this function $f$ we get $h_{\ell}(\a)=\a^{-1}$.

Having  possibility to choose several functions in construction of
GDPR estimators, natural question is what properties of these
functions ensures better results in estimating $\a$. Comparing two
functions $f_1(x)=x$ and $f_{\ell}(x)=-\log x$ we see that
corresponding functions $h_1(\a)=\a/(1+\a)$ and
$h_{\ell}(\a)=\a^{-1}$ have quite different ranges: the first one
has a small range - interval $(0,1)$, while the second one as a
range have all half line $(0, \infty)$. Moreover, this fact results
in different behavior of derivatives of inverse functions
\begin{eqnarray*}
\frac{{\rm d} }{{\rm d} p}h_1^{-1}(p) &=&\frac{{\rm d}}{{\rm d} p}\left( \frac{p}{1-p}\right)=\frac{1}{(1-p)^2}=(\alpha+1)^2, \\
\frac{{\rm d} }{{\rm d} \gamma} h_{\ell}^{-1}(\gamma)&=&\frac{{\rm
d} }{{\rm d} p}
\left(\frac{1}{\gamma}\right)=-\frac{1}{\gamma^2}-=\alpha^2.
\end{eqnarray*}
For small values of $\a$ (this corresponds to small values of $p$
and large values of $\gamma$) the derivative of the first function
is almost one, while for the second function it tends to zero as
$\gamma^{-2}$. This means that even big changes in the value of
estimated quantity $\gamma$ results only in small changes of
estimated value of $\a$. For large values of $\a$ (as $p\to 1$ or
$\gamma \to 0$) behavior of both derivatives is the same, but,
evidently, large values of $\a$ are not so interesting in the
problem of tail index estimation. These considerations explain why
Qi estimator (\ref{Qi-s}) with $s=1$ (or, in other words, GDPR
estimator with the function $f_{\ell}$) performs better than DPR
estimator (with the function $f_1$) and also suggest one more
modification of DPR estimators. Namely, if we take the same function
$f_r$, but with negative parameter $r$ (to stress this we shall
write $f_{-r}, \ r>0$), there will appear restriction $\a>r$, but
now the range of the function $h_{-r}(\a)=\a(\a-r)^{-1}$ is infinite
interval $(1, \infty)$ and the behavior of the derivative of inverse
function is very similar to that of $h_{\ell}^{-1}$. We are able to
show that in the case of function $f_{-r}$ it is possible to find
optimal choice of $r$ and GDPR estimator with this optimal $r$ is
comparable with estimator (\ref{Qi-s}) with $s=1$  in a sense that
for some values of $\a, \ \b$ one estimator has smaller asymptotic
MSE, for other values - dominates another one.
 Investigation of all these modifications were carried while the
 first version of the paper was in the process of refereeing and the
 results with proofs are collected in a forthcoming paper
 \textsc{Paulauskas} and \textsc{Vai\v ciulis} (2010).

One more remark concerning Theorem \ref{thm2} is appropriate here.
In \textsc{Qi} (2010) it is mentioned that using the same method of
the proof of asymptotic normality for estimator (\ref{Qi-s})(that
is, using the relation between ordered statistics and exponential
distributions) it is possible to prove (\ref{clt_eq}). Our proof of
(\ref{clt_eq}) essentially differs since it does not use exponential
distributions and the main tool in the proof is formula (\ref{new}).
It is worth to mention that one can get the results for estimator
(\ref{Qi-s}) with $s=1$ by using this formula. This is demonstrated
in the above mentioned forthcoming paper.

The paper contains two more sections. In Section 2 we prove Theorem
\ref{thm2} and in Section 3 there are results on comparison of
estimators.

\section{Proof of Theorem \ref{thm2}}

{\it Proof of (\ref{bias_ratio_asymp}).} Relation
(\ref{bias_ratio_asymp}) gives the exact asymptotic of the bias $\E
\hat{p} - p$. Generally, the exact asymptotic of the bias of a tail
index estimator is rather difficult problem. The advantage of our
estimator $\hat{p}$ is a relative simplicity of the  proof of
(\ref{bias_ratio_asymp}).
 We do not use asymptotic for the inverse function for
$1-F(x)$  as in \textsc{Paulauskas} (2003), but rather simple form of the
expectation
\begin{eqnarray}
\E \hat{p} &=& 1- m \int_0^{\infty} F^{m-1}(x_2) \left\{
\int_{x_2}^{\infty} \frac{\dd F(x_1)}{x_1} \right\} \dd x_2,
\label{new}
\end{eqnarray}
which will be proved below. We  truncate the outer integral at the
level
\begin{equation} \label{a_m}
 a_m= \kappa
m^{1/\alpha} \left( \ln m \right)^{-1/\alpha},
\end{equation}
where $0< \kappa<(C_1/\zeta)^{1/\alpha}$ and denote
\begin{eqnarray*}
K_{m,1} &=& m \int_0^{a_m} F^{m-1}(x_2) \left\{ \int_{x_2}^{\infty} \frac{\dd F(x_1)}{x_1} \right\} \dd x_2, \\
K_{m,2} &=& m \int_{a_m}^{\infty} F^{m-1}(x_2) \left\{ \int_{x_2}^{\infty} \frac{\dd F(x_1)}{x_1} \right\} \dd x_2, \\
\end{eqnarray*}
Now, (\ref{bias_ratio_asymp}) follows immediately from the following two relations
\begin{eqnarray} \label{rel_01}
K_{m,1} &=& o\left(m^{-\zeta}\right), \\
\label{rel_02} 1-p-K_{m,2} &\sim& \chi m^{-\zeta}, \quad m \to
\infty.
\end{eqnarray}
To prove (\ref{rel_01}), split $K_{m,1}$ into two parts:
$K_{m,1}=K'_{m,1}+K''_{m,1}$, where
$$
K'_{m,1} = m \int_0^{a_m} F^{m-1}(x_2) \left\{ \int_{x_2}^{a_m}
\frac{\dd F(x_1)}{x_1} \right\} \dd x_2.
$$
By the change of integration order to get
\begin{eqnarray*}
K'_{m,1} &=& m \int_0^{a_m}  \frac{\dd F(x_1)}{x_1} \left\{ \int_{0}^{x_1} F^{m-1}(x_2) \right\} \dd x_2 \\
&\le& m \int_0^{a_m}  F^{m-1}(x_1) \dd F(x_1)  =F^m(a_m).
\end{eqnarray*}
An assumption (\ref{secorder}) and a simple inequality $\ln(1-x) \le
-x$, $0 \le x<1$ yield
$$
F^m(a_m) \le C\left(1-C_1 a_m^{-\alpha}\right)^m \le C {\rm e}^{m
\ln (1-C_1 a_m^{-\alpha})}  \le C m^{-C_1 \kappa^{-\alpha}}
$$
 for sufficiently large $m$, hence, taking into account (\ref{a_m}), we get $K'_{m,1} =o\left( m^{-\zeta}\right)$.
  Relations
$\int_0^{a_m} F^{m-1}(x) \dd x =O\left(a_m F^{m-1}(a_m) \right)$ and
$\int_{a_m}^{\infty} x^{-1} \dd F(x)=O\left(a_m^{-\alpha-1} \right)$
prove $$K''_{m,1}=O\left( m a_m^{-\alpha} F^{m-1}(a_m)
\right)=o\left( m^{-\zeta}\right),$$ and we have (\ref{rel_01}).

 Now let us prove (\ref{rel_02}). Integrating by parts the inner integral
we get
\begin{eqnarray*}
K_{m,2} = m \int_{a_m}^{\infty} F^{m-1}(x) \left\{
\frac{C_1 \alpha}{\alpha+1} x^{-\alpha-1}+  \frac{ C_2 \beta}{\beta+1}
x^{-\beta-1}+o\left( x^{-\beta-1}\right) \right\} \dd x.
\end{eqnarray*}
Denote $\tilde{F}(x)=1-C_1 x^{-\alpha} -C_2 x^{-\beta}$, then one
can write
$$
K_{m,2}=K'_{m,2}+K''_{m,2}+R_{m},
$$
where
\begin{eqnarray}
K'_{m,2} &=& \frac{m}{\alpha+1} \int_{a_m}^{\infty} \tilde{F}^{m-1}(x)   \dd \tilde{F}(x), \nonumber \\
\label{K''_{m,2}} K''_{m,2} &=&  C_2 \beta m \left(
\frac{1}{\beta+1} - \frac{1}{\alpha+1}\right)  \int_{a_m}^{\infty}
x^{-\beta-1} \tilde{F}^{m-1}(x)  \dd x,
\end{eqnarray}
and \begin{equation}\label{rem1}R_m=R_{m,1}+R_{m,2}\end{equation}
with
\begin{eqnarray*}
R_{m,1} &=& m \int_{a_m}^{\infty} \left( F^{m-1}(x)-
\tilde{F}^{m-1}(x)\right) \left\{ \frac{C_1\alpha}{\alpha+1}
x^{-\alpha-1}+ \frac{C_2 \beta}{\beta+1}
x^{-\beta-1} \right\} \dd x, \\
R_{m,2} &=& m \int_{a_m}^{\infty}  F^{m-1}(x) o\left(
x^{-\beta-1}\right) \dd x.
\end{eqnarray*}
Integration gives $1-p-K'_{m,2}=O\left(
\tilde{F}^m(a_m)\right)=o\left( m^{-\zeta}\right)$. We claim that
\begin{equation}\label{cl0}
    K''_{m,2} \sim - \chi m^{-\zeta}, \quad
    (m \to \infty).
\end{equation}
Simple calculations show that (\ref{cl0}) can be derived  from the
following two relations:
\begin{eqnarray}
\label{cl1} \int_{a_m}^{\infty} \exp\{-C_1 (m-1) x^{-\alpha}\}
x^{-\beta-1} \dd x \sim \frac{\Gamma(\beta/\alpha)}{\alpha C_1^{\beta/\alpha}} m^{-\beta/\alpha}, \\
 \label{cl2}
\int_{a_m}^{\infty} \left\{ \tilde{F}^{m-1}(x) -\exp\{-C_1 (m-1)
x^{-\alpha}\}\right\} x^{-\beta-1} \dd x = o\left(
m^{-\beta/\alpha}\right).
\end{eqnarray}
Making a change of  variables $y=C_1 (m-1) x^{-\alpha}$ one has
\begin{eqnarray*}
&&    \int_{a_m}^{\infty} \exp\{-C_1 (m-1) x^{-\alpha}\}
x^{-\beta-1} \dd x \\
&& \quad = \frac{1}{\alpha C_1^{\beta/\alpha}} (m-1)^{-\beta/\alpha}
\left(\Gamma(\beta/\alpha)- \Gamma(\beta/\alpha, C_1(m-1)
a_m^{-\alpha}) \right),
\end{eqnarray*}
where $\Gamma(\cdot)$ is a standard  gamma function and
$$
\Gamma(\alpha, x)=\int_x^\infty t^{\alpha-1}e^{-t}\dd t
$$
is an upper incomplete gamma function. Keeping in mind that
$a_m^{-\alpha} m \to \infty$, we have
\begin{eqnarray*}
\Gamma(\beta/\alpha, C_1(m-1) a_m^{-\alpha}) \le \frac{1}{C_1(m-1)
a_m^{-\alpha}} \Gamma\left(1+\frac{\beta}{\alpha}\right) \to 0,
\quad (m \to \infty).
\end{eqnarray*}
This ends the  proof of (\ref{cl1}).

 To prove (\ref{cl2}), consider the difference
$\Delta_{m}(x):=\tilde{F}^{m}(x) -\exp\{-C_1 m x^{-\alpha}\}$. We
assume that  $m$ is large enough that inequalities
\begin{equation}\label{hk1}
0<C_1 a_m^{-\alpha}+C_2 a_m^{-\beta}<1/2, \quad 1+\frac{C_2 \beta}{C_1 \alpha} a_m^{-\beta+\alpha}>0
\end{equation}
are satisfied. We recall that only $C_1$ is supposed positive, while
$C_2$ may be negative. The second inequality in (\ref{hk1}) ensures
monotonous decay of a function $1-\tilde{F}(x)$, $x\ge a_m$ and
implies $0<1-\tilde{F}(x)<1/2$ for $x\ge a_m$. Then we can write :
$$
\ln \tilde{F}(x) = -C_1 x^{-\alpha} -C_2 x^{-\beta} + r(x),
$$
where $|r(x)| \le 2\left(C_1+|C_2| \right)^2 x^{-2 \alpha}$. Using
this relation, for sufficiently large $m$, we have
\begin{eqnarray*}
\left| \Delta_{m}(x) \right| &=& \exp\{-C_1 m x^{-\alpha}\}\left|\exp\left\{ C_1 m x^{-\alpha} +m \ln \tilde{F}(x)\right\} -1\right| \\
&=& \exp\{-C_1 m x^{-\alpha}\} \left|\exp\left\{ -m\left(C_2 x^{-\beta} - r(x)\right)\right\} -1\right| \\
&\le&  C m x^{-(2\alpha \wedge \beta)} \exp\{-C_1 m x^{-\alpha}\}.
\end{eqnarray*}
 Consequently, left hand side of
(\ref{cl2}) does not exceed
\begin{equation} \label{cl2a}
m \int_{a_m}^{\infty}  \exp\{-C_1 (m-1) x^{-\alpha}\} x^{-(2\alpha
\wedge \beta)-\beta-1} \dd x= o\left(m^{1-\frac{(2\alpha \wedge
\beta)+\beta}{\alpha}} \right).
\end{equation}
If $2\alpha \ge \beta$, then r.h.s. of (\ref{cl2a}) is
$o\left(m^{1-2\beta/\alpha}\right)=o\left(m^{-\beta/\alpha}\right)
$, while, in the case $2\alpha < \beta$, we have
$o\left(m^{1-(2\alpha \wedge \beta)/\alpha-\beta/\alpha} \right) =
o\left(m^{-\beta/\alpha-1}\right)$. Thus, (\ref{cl2}) and,
consequently, (\ref{cl0}), are proved.

To finish the proof of (\ref{rel_02}), it remains to prove that the
remainder term $R_m$ from (\ref{rem1}) is negligible, that is,
$R_{m,i}=o\left(m^{-\zeta} \right), \ i=1,2$. We have
$$
\left| R_{m,1}\right| \le  C m \int_{a_m}^{\infty} \left| F^{m-1}(x)- \tilde{F}^{m-1}(x)\right| x^{-\alpha-1}\dd x.
$$
 The remainder term in (\ref{secorder}) denote
by $h(x)$, that is,
$h(x)=1-F(x)-C_1x^{-\alpha}-C_2x^{-\beta}=x^{-\beta}h_1(x)$ where
$h_1(x)=o(1)$, as $x \to \infty$. Let us rewrite the difference in
the integrand as follows
$$
F^{m-1}(x)- \tilde{F}^{m-1}(x) = \tilde{F}^{m-1}(x)
\left(\exp\left\{(m-1) \ln \left(1-\frac{h(x)}{\tilde{F}(x)} \right)
\right\} -1\right).
$$
One can assume that for $x \ge a_m$  inequality $\left|h(x)
/\tilde{F}(x)\right|<1/2$ is satisfied, thus from the Taylor
expansion of $\ln \left(1-h(x)/\tilde{F}(x) \right)$ it follows that
there exist a constant $C>0$ such that $|\ln
\left(1-h(x)/\tilde{F}(x) \right)| \le C |h(x)|/\tilde{F}(x)$. Since
for $x \ge a_m$ the product $(m-1) x^{-\beta}$ tends to zero as $m
\to \infty$, we can assume that inequality
$$
(m-1) \left|\ln \left(1+\frac{h(x)}{\tilde{F}(x)} \right)\right| <1
$$
holds. Then, by applying inequality $\left| {\rm e}^x-1\right| \le C|x|$, $0<|x|<1$, we get
\begin{eqnarray*}
\left|F^{m-1}(x)- \tilde{F}^{m-1}(x)\right| &\le& (m-1) \tilde{F}^{m-1}(x)  \left|\ln \left(1-\frac{h(x)}{\tilde{F}(x)} \right) \right| \\
&\le&C (m-1) \tilde{F}^{m-1}(x)  \left|h(x)\right|.
\end{eqnarray*}
Applying the last inequality we have
\begin{eqnarray*}
R_{m,1} &\le& C m \int_{a_m}^{\infty} \tilde{F}^{m-1}(x)  \left|h(x)\right| x^{-\alpha-1}\dd x \\
&\le&  C m a_m^{-\alpha} \sup_{x\ge a_m} \left|
h_1(x)\right|\int_{a_m}^{\infty} x^{-\beta-1} \tilde{F}^{m-1}(x) \dd
x.
\end{eqnarray*}
Taking into account  (\ref{cl1})-(\ref{cl2}) we obtain  $R_{m,1} \le
C a_m^{-\alpha} m^{-\zeta} \sup_{x\ge a_m}\left|
h_1(x)\right|=o\left(m^{-\zeta}\right)$. In a similar way we get
$$
\left|R_{m,2}\right| \le C m^{-\zeta} \sup_{x \ge a_m} |h(x)| =o\left(m^{-\zeta}\right)
$$
and the proof of (\ref{rel_02}) is completed.

To complete the proof of (\ref{bias_ratio_asymp}) it remains to
prove (\ref{new}). The random variables $\kappa_{n,1},\dots,
\kappa_{n,n}$, defined in (\ref{definitions}), are i.i.d.. Therefore
$\E\hat p = \E \kappa_{n,1}$.
Now it is clear that
\begin{eqnarray}
\E\hat p  &=& m! \int_{0}^{\infty} \frac{\dd F(x_1)}{x_1} \int_{0}^{x_1} x_2 \dd F(x_{2}) \int_{0}^{x_2} \dd F(x_{2}) \dots \int_{0}^{x_{m-1}} \dd F(x_{m}) \nonumber \\
&=&(m-1)m \int_{0}^{\infty} \left\{ \int_{0}^{x_1} x_2 F^{m-2}(x_2) \dd F(x_{2}) \right\}\frac{\dd F(x_1)}{x_1} .
\nonumber
\end{eqnarray}
Integrating by parts the inner integral we get
\begin{eqnarray}
\E\hat p &=& m \int_{0}^{\infty} \left\{ x_1 F^{m-1}(x_1) - \int_{0}^{x_1} F^{m-1}(x_2) \dd x_{2} \right\}\frac{\dd F(x_1)}{x_1} \nonumber \\
&=& 1- m \int_{0}^{\infty} \left\{\int_{0}^{x_1} F^{m-1}(x_2) \dd x_{2} \right\}\frac{\dd F(x_1)}{x_1}. \nonumber
\end{eqnarray}
It remains  to change order of integration to conclude
the proof of (\ref{new}).

\bigskip

\noindent {\it Proof of (\ref{clt_eq}).} Since the proofs of (\ref{clt_eq}) and (\ref{clt_paul})
are essentially the same, we shall give main steps only. In view of decomposition
$$
\sqrt{n} \left( \hat{p}-p\right) = \frac{1}{\sqrt{n}}\sum_{i=1}^n
\left (\kappa_{n, i} - \E \kappa_{n,i} \right) + \sqrt{n} \gamma_m,
$$
the assertion (\ref{clt_eq}) follows from
\begin{eqnarray}
\label{step1} \frac{1}{\sqrt{n}}\sum_{i=1}^n \left (\kappa_{n, i} -
\E \kappa_{n,i} \right) \buildrel \rm D \over \longrightarrow_{N\to
\infty}
N(0,\sigma^2), \\
\label{step2} \sqrt{n} \gamma_m \to \mu, \quad (m=m_{opt}, \ N \to
\infty).
\end{eqnarray}
To prove (\ref{step1}) check Lyapunov condition  for i.i.d. random variables forming
triangular array, while relation
(\ref{step2}) one can verify by using (\ref{bias_ratio_asymp}) and
(\ref{opt4}) with $m=m_{opt}$.

\bigskip

\noindent {\it Proof of (\ref{opt5}).} From (\ref{clt_eq}) we know
that the variance $\E\left( \hat{p} - \E \hat{p}\right)^2$
asymptotically equals $\sigma^2 m/N$. Taking the main term in the
asymptotic relation (\ref{bias_ratio_asymp}) we get that the
asymptotic mean squared error of $\hat{p}$ equals $\chi^2 m^{-2
\zeta}+\sigma^2 m/N$. Since this is a very simple function with
respect to $m$ it is easy to see that a solution of a minimization
problem
\begin{equation} \label{minimiz}
\inf_{2 \le m \le N}\left\{\chi^2 m^{-2 \zeta}+\sigma^2 m/N \right\}
\end{equation}
is given by (\ref{opt4}). Here it is necessary to note that we
require that $N$ is sufficiently large, since for  values of $\zeta$
close to $0$ the solution of minimization problem (for a given $N$)
may be smaller than $2$. Also one can note that instead of taking
main term from (\ref{bias_ratio_asymp}) we can take the sequence
$\gamma_m$ and apply Lemma  2.8 in  \textsc{Dekkers} and \textsc{de Haan} (1993), this will give
the same result. Having (\ref{opt4})one can easily get (\ref{opt5}).
The theorem is proved.

\section{Comparison of estimators }

In this section we compare the  tail index estimator $\hat p$ with
the estimators $\gamma_{N,k}^{(j)}$, $j=1,2,3,4$, using the same
method as in \textsc{de Haan} and  \textsc{Peng} (1998). May be here it is worth to mention that in
a recent paper \textsc{Nematollahi} and \textsc{Tafakori} (2007)
there was proposed another approach to
compare tail index estimators, but this method of comparison is well
adapted to a specific estimator introduced by \textsc{fan} in  \textsc{fan} (2004). We
recall that $\hat p$ estimates the quantity  $p=\alpha/(\alpha+1)$,
while the four above mentioned estimators estimate
$\gamma=1/\alpha$, that is, different function of the unknown
parameter $\a$. Therefore the first step in comparison is to
transfer the estimators to the same function, and we had chosen  to
compare the estimator $\hat{p}$ with estimators
$p_{N,k}^{(j)}=1/(1+\gamma_{N,k}^{(j)})$. We need the following
simple statement.

\begin{lem} \label{lema1} Suppose (\ref{secorder}) holds.
Let $k^{(j)}=k^{(j)}(N)$, $j=1,2,3,4$ be a sequences of integers with
\begin{equation}\label{k-req}
k^{(j)}(N) \to \infty \quad {\rm and} \quad k^{(j)}(N)/N \to 0, \quad (N \to \infty)
\end{equation}
and let estimators $\gamma_{N,k}^{(j)}$
are asymptotically normal, i.e. there exist constants $b_j \in \mathbb{R}$ and $\sigma_j >0$ such that
\begin{equation}\label{CLT_gamma}
    \sqrt{k^{(j)}}\left( \gamma_{N,k}^{(j)} -\gamma\right)\buildrel \rm D \over \longrightarrow_{N\to \infty} \mathcal{N}(b_j, \sigma^2_j).
\end{equation}
Then
\begin{equation}\label{CLT_p}
    \sqrt{k^{(j)}}\left( p_{N,k}^{(j)} -p\right)\buildrel \rm D \over \longrightarrow_{N\to \infty}
\mathcal{N}\left(-\frac{b_j}{(1+\gamma)^2}, \frac{\sigma^2_j}{(1+\gamma)^4}\right).
\end{equation}
\end{lem}
{\it Proof.} We use the
obvious identity
$$
\sqrt{k^{(j)}}\left( p_{N,k}^{(j)} -p\right) =
\frac{\sqrt{k^{(j)}}\left( \gamma-\gamma_{N,k}^{(j)}
\right)}{(1+\gamma)(1+\gamma_{N,k}^{(j)})}.
$$
Now  relation  (\ref{CLT_p}) follows from this identity, Theorem 4.4
in \textsc{Billingsley} (1968), relation (\ref{CLT_gamma}) and the relation
$\gamma_{N,k}^{(j)}\buildrel \rm P \over \rightarrow \gamma$, ($N \to \infty$). The lemma is proved.

\bigskip

De Haan and Peng proved (see Thm.2 in \textsc{de Haan} and  \textsc{Peng} (1998)) that the
asymptotic second moment of $\gamma_{N,k}^{(j)} -\gamma$ is minimal
and equals $\left( k^{(j)} \right)^{-1} \sigma^2_k (1+2\zeta)/(2
\zeta)$, if
 $k^{(j)}$ satisfies the  relation
\begin{equation}\label{extr_k}
    \lim_{N \to \infty} k^{(j)} A^2(N/k^{(j)}) = \frac{\sigma_j^2}{2 \zeta D_k^2
    }, \quad j=1, 2, 3, 4,
\end{equation}
where
\begin{eqnarray*}
D_1 &=& \frac{1}{1+\zeta}, \ D_3 = \frac{1}{1+\zeta} -\frac{\alpha \zeta}{(1+\zeta)^2},\\
D_2 &=& \frac{1}{(2^{1/\alpha} -1)\ln 2} \frac{1-2^{-\zeta}}{\zeta}
\left(2^{(1/\alpha)-\zeta} -1\right), \ D_4 = \frac{1}{(1+\zeta)^2},
\end{eqnarray*}
and
\begin{eqnarray*}
\sigma_1^2 =  \frac{1}{\alpha^2}, \
\sigma_2^2 =  \frac{1+2^{2/\alpha+1}}{\alpha^2 (2^{1/\alpha}-1)^2 \ln^2 2} , \
\sigma_3^2 =  \frac{1+\alpha^2}{\alpha^2}, \
\sigma_4^2 =  \frac{2}{\alpha^2}
\end{eqnarray*}
are limit variances in (\ref{CLT_gamma}). The function $A(t)$ in
(\ref{extr_k}) has the asymptotic
$$
A(t) \sim -\frac{\zeta}{\alpha} \frac{C_2}{C_1^{\beta/\alpha}}
t^{-\zeta}, \quad t \to \infty.
$$
We recall  that the function $A(t)$ was introduced in
(\ref{gensecor}).

Denote by $k^{(j)}_{opt}$ a sequence $k^{(j)}$, satisfying (\ref{extr_k}). From (\ref{extr_k}), we have
\begin{eqnarray*}
k^{(1)}_{opt}(N) &\sim& \left( \frac{(1+\zeta)^2}{2 \zeta^3} \frac{\left(C_1 \right)^{2\beta/\alpha}}{\left(C_2 \right)^{2}}\right)^{1/(1+2\zeta)} N^{2\zeta/(1+2\zeta)},\\
k^{(2)}_{opt}(N) &\sim& \left( \frac{1+2^{2/\alpha+1}}{2 \zeta (1-2^{-\zeta})^2 (2^{1/\alpha-\zeta}-1)^2} \frac{\left(C_1 \right)^{2\beta/\alpha}}{\left(C_2 \right)^{2}}\right)^{1/(1+2\zeta)} N^{2\zeta/(1+2\zeta)},\\
k^{(3)}_{opt}(N) &\sim& \left( \frac{(1+\zeta)^4 (1+\alpha^2)}{2 \zeta^3 (1+\zeta -\zeta \alpha)^2} \frac{\left(C_1 \right)^{2\beta/\alpha}}{\left(C_2 \right)^{2}}\right)^{1/(1+2\zeta)} N^{2\zeta/(1+2\zeta)},\\
k^{(4)}_{opt}(N) &\sim& \left( \frac{(1+\zeta)^4}{ \zeta^3} \frac{\left(C_1 \right)^{2\beta/\alpha}}{\left(C_2 \right)^{2}}\right)^{1/(1+2\zeta)} N^{2\zeta/(1+2\zeta)}.
\end{eqnarray*}
Under normalization $k^{(j)}_{opt}(N)$ instead of $k^{(j)}(N)$ in
(\ref{CLT_gamma}) the limit random variable has a mean
$$
\frac{\sigma_k}{\sqrt{2 \zeta}}{\rm sgn}\left(D_k \lim_{N \to
\infty} \sqrt{k^{(j)}_{opt}(N)}
A\left(N/k^{(4)}_{opt}(N)\right)\right).
$$
Moreover, Lemma  \ref{lema1} imply
\begin{equation}\label{opt_mse_p}
    \E\left(p_{N,k}^{(j)}-p \right)^2 \sim \left(\frac{\alpha}{\alpha+1}\right)^4  \frac{2\beta-\alpha}{2(\beta-\alpha)} \frac{\sigma^2_j}{k^{(j)}_{opt}(N)}, \quad (N \to \infty).
\end{equation}

Now it is possible to compare the estimator $\hat p$ with the
estimators $p_{N,k}^{(j)}$  as in was done in \textsc{de Haan} and  \textsc{Peng} (1998), i.e., by
calculating a limit of the ratio of minimal mean squared errors:
$$
RMMSE(j)=\lim_{N \to \infty} \frac{\E\left(\hat{p}-p \right)^2}{\E\left(p_{N,k}^{(j)}-p \right)^2}.
$$
From (\ref{opt5}) and (\ref{opt_mse_p}) we have the following
results:
\begin{eqnarray*}
RMMSE(1) &=& \left( \eta(\alpha, \beta) \Gamma^2(2+\zeta)   \right)^{1/(1+2\zeta)}, \\
RMMSE(2) &=& (2^{(1/\alpha)}-1)^2 \ln^2(2) \bigg(
\eta(\alpha, \beta) \times \\
&& \times \frac{\zeta^2 \Gamma^2(1+\zeta) }{(1-2^{-\zeta})^2(2^{(2/\alpha)+1}+1)^{2\zeta} (2^{1/\alpha-\zeta}-1)^2}  \bigg)^{1/(1+2\zeta)},\\
RMMSE(3) &=& \left( \eta(\alpha, \beta) \frac{(1+\zeta)^{2} \Gamma^2(2+\zeta)}{(1+\alpha^2)^{2\zeta}(1+\zeta-\alpha \zeta)^2}
 \right)^{1/(1+2\zeta)}, \\
RMMSE(4) &=&
\left( \eta(\alpha, \beta) \frac{(1+\zeta)^{2} \Gamma^2(2+\zeta)}{2^{2\zeta}}  \right)^{1/(1+2\zeta)},
\end{eqnarray*}
where
$$
\eta(\alpha, \beta)=\left(\frac{\beta(\alpha+1)}{\alpha(\beta+1)}\right)^2
\left(\frac{(\alpha+1)^2}{\alpha(\alpha+2)}\right)^{2\zeta}.
$$
It is easy to conclude that $RMMSE(1)>1$ for all $0<\alpha<\beta$
i.e., Hill estimator $p_{N,k}^{(1)}$  dominates  estimator
$\hat{p}$. Due to the inequality $(\alpha+1)^6 - 4
\alpha^3(\alpha+2)>0$, for $\alpha>0$ (it follows from the binomial
formula), the same conclusion is valid for de Vries estimator
$p_{N,k}^{(4)}$ .

Comparison of estimators $\hat{p}$, $p_{N,k}^{(2)}$ and
$p_{N,k}^{(3)}$ is shown in Figures 1a-1c. $\alpha$ values are on
the horizontal axis, while vertical axis labels $\beta$ values. In
all three figures the area $\{(\alpha, \beta): \ 0<\beta <\alpha\}$
(those values of parameters that are not considered) is left as
white. In Figure 1a the area $\{ (\alpha, \beta): \ RMMSE(2)>1\}$ is
in black and $\{ (\alpha, \beta): \ RMMSE(2)<1\}$ is in dark grey.
Similarly, Figure 1b presents comparison of the estimators $\hat{p}$
and $p_{N,k}^{(3)}$: as in Figure 1a, the area $\{ (\alpha, \beta):
\ RMMSE(3)<1\}$ is in dark grey and the area $\{ (\alpha, \beta): \
RMMSE(3)>1\}$ - in light grey. Finally, Figure 1c gives areas of
domination estimators $\hat{p}$ (dark grey), $p_{N,k}^{(2)}$
(black), and $p_{N,k}^{(3)}$ (light grey).

\begin{figure}[ht!] \label{alpha}
\begin{center}
\includegraphics[width=0.30\textwidth]{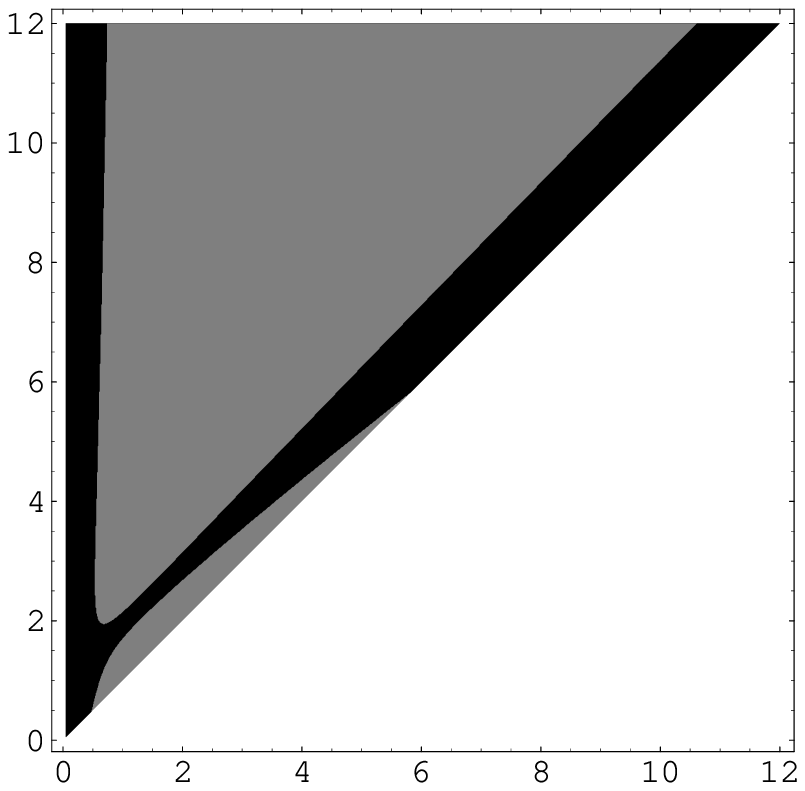} \quad \includegraphics[width=0.30\textwidth]{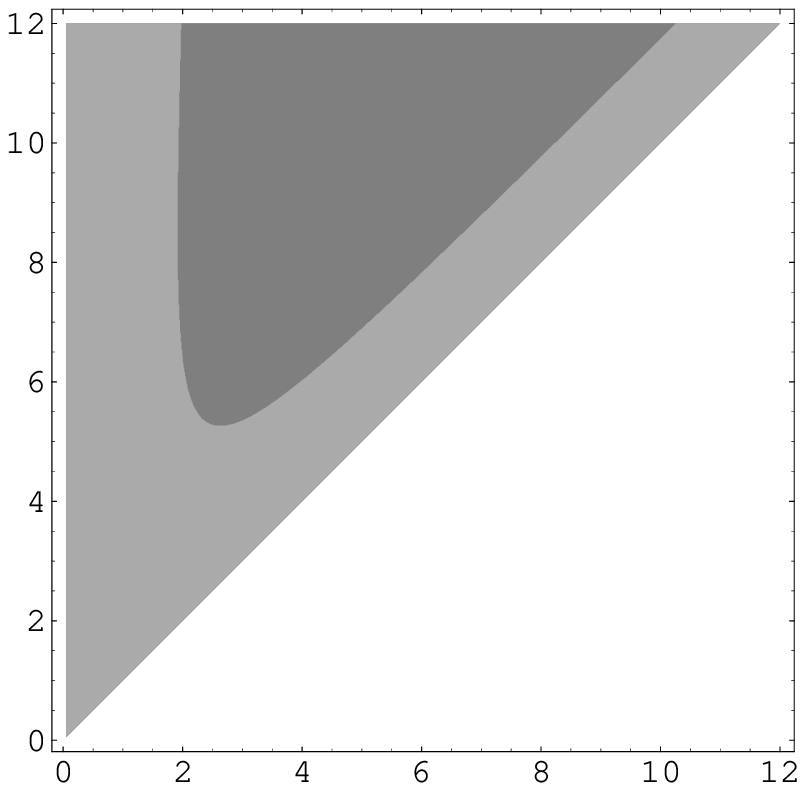} \quad \includegraphics[width=0.30\textwidth]{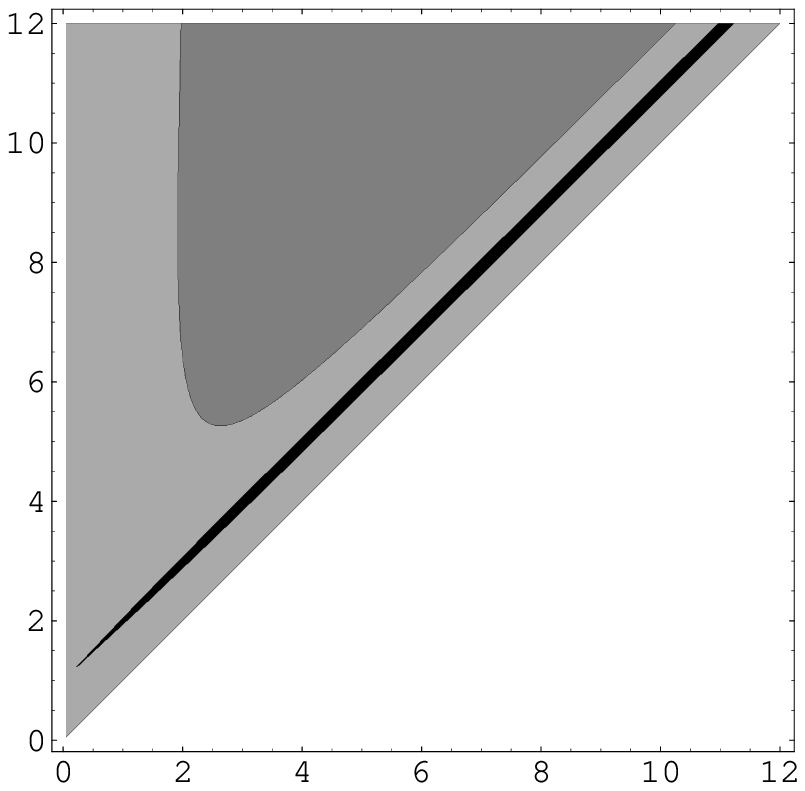}
\end{center}
\caption{1a, 1b and 1c figures}
\end{figure}

As it was mentioned  in the Introduction, in \textsc{de Haan} and  \textsc{Peng} (1998) the
comparison of the estimators $\gamma_{N,k}^{(j)}$, $j=1,2,3,4$ was
performed with respect to the parameters $(\gamma, \rho)$. For the
sake of completeness we include analogous of the Figures 1a-1c in
the plane $(\rho, \gamma)$ also. In the  Figures 2a-2c the
horizontal axis labels $\gamma$ values, while vertical axis labels
$\rho$ values. As in Figures 1a-1c, the area where estimator
$\hat{p}$ has an asymptotic mean squared error smaller than the
other estimator(s) is in dark grey. A black and light grey colors
mark the areas of domination of estimators $p_{N,k}^{(2)}$ or
$p_{N,k}^{(3)}$, respectively.

\begin{figure}[ht!] \label{alpha1}
\begin{center}
\includegraphics[width=0.30\textwidth]{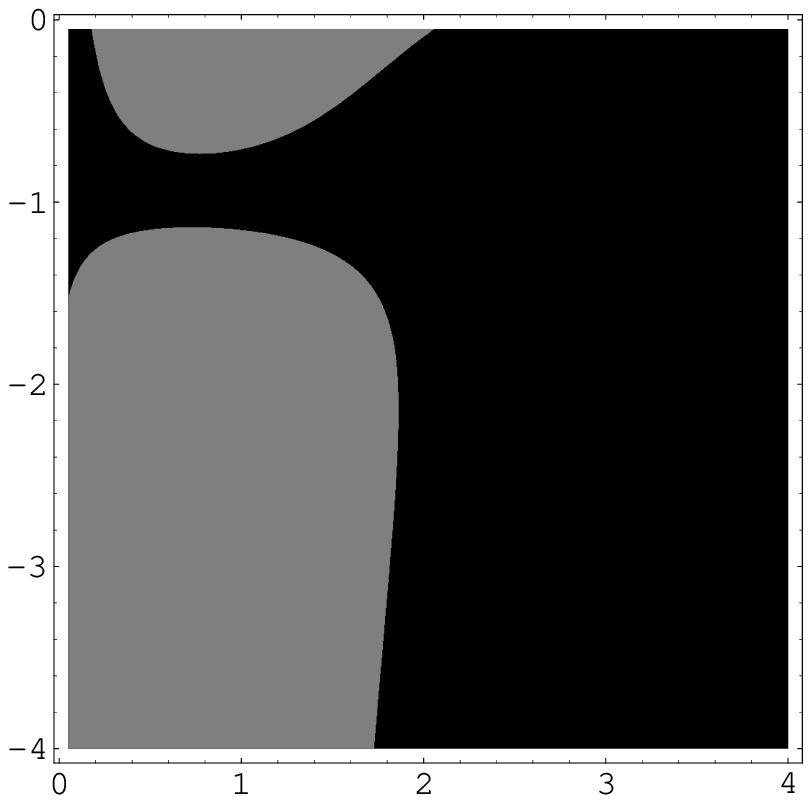} \quad \includegraphics[width=0.30\textwidth]{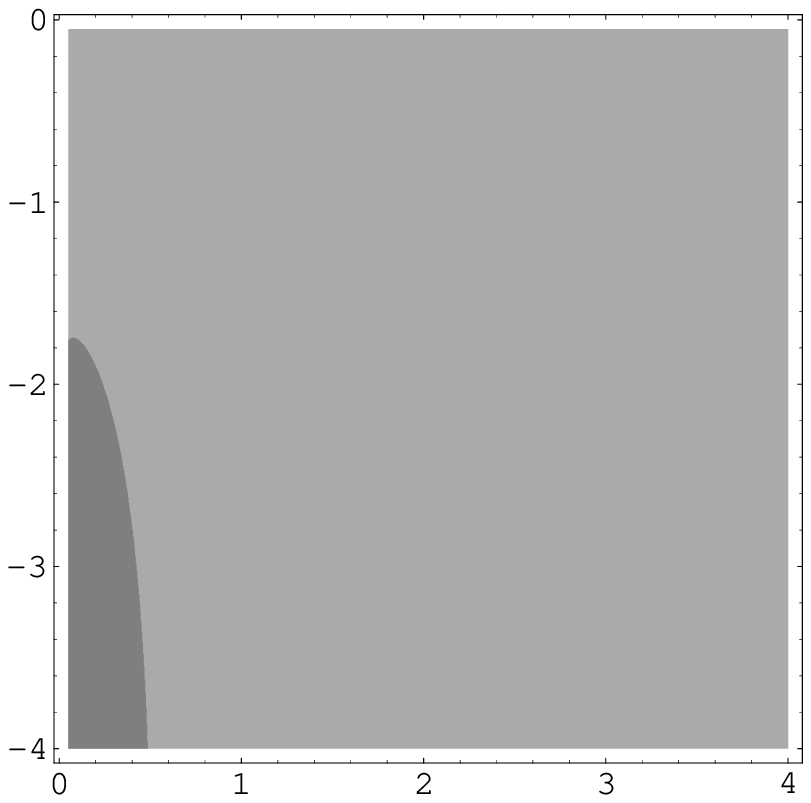} \quad \includegraphics[width=0.30\textwidth]{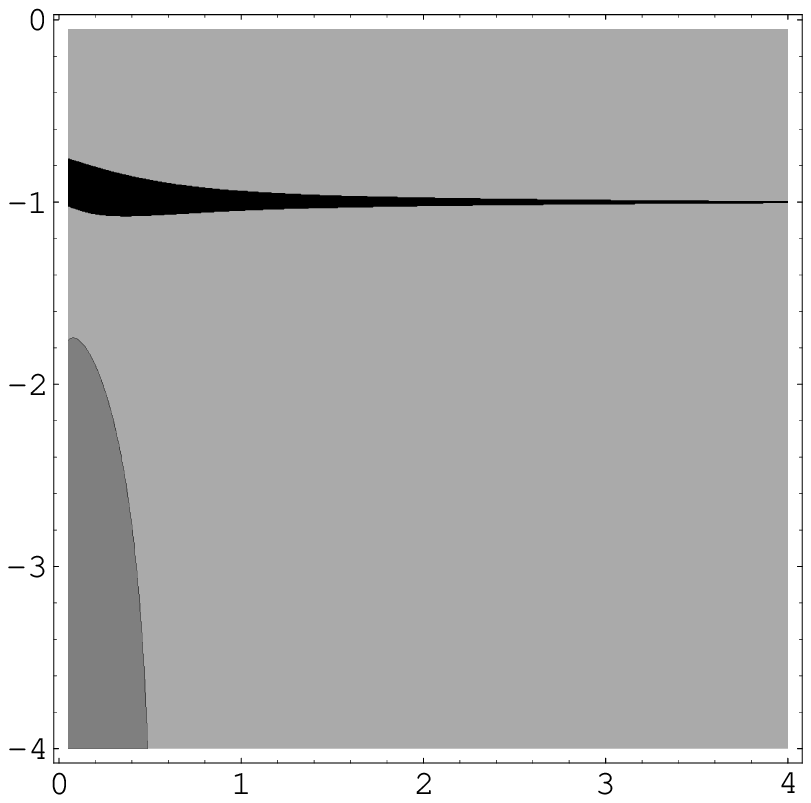}
\end{center}
\caption{2a, 2b and 2c figures}
\end{figure}

\newpage

\section*{References}

\begin{description}

\itemsep -.04cm

\item{\sc Billingsley, P.} (1968), {\em Convergence of Probability Measures}, Wiley, New York.

\item{\sc Davydov, Yu. and V. Paulauskas} (1999), On the estimation of the parameters of multivariate stable distributions,
{\em Acta Appl. Math.}, 58, 107--124.

\item{\sc Davydov, Yu.,   V. Paulauskas and A. Ra\v ckauskas} (2000), More on $P$-stable convex sets in Banach spaces,
{\em J. Theoret. Probab.}, 13, 39--64.

\item{\sc De Haan, L. and  L. Peng} (1998), Comparison of tail index estimators,
{\em Statist. Neerlandica}, 52, 60--70.

\item{\sc Dekkers, A.L.M. and L. de Haan} (1993), Optimal choice of sample fraction in extreme-value estimation,
{\em J. Multivariate Anal.}, 47, 173--195.

\item{\sc Dekkers, A.L.M.,  J.H.J. Einmahl and   L. de Haan} (1989), A moment estimator for the index of an extreme-value distribution,
{\em Ann. Statist.}, 17, 1833-1855.

\item{\sc Fan, Zh.} (2004), Estimation problems for distributions with heavy tails,
{\em J. Statist. Plann. Inference}, 123, 13--40.

\item{\sc Fraga Alves, M.I.} (2001), A location invariant {H}ill-type estimator,
{\em Extremes}, 4, 199--217.

\item{\sc  Fraga Alves, M.I.,  L. de Haan and T. Lin} (2006), Third order extended regular variation,
{\em Publications de lInstitut Math´ematique}, 80, 109--120.

\item{\sc Fraga Alves, M.I., M.I. Gomes and  L. de Haan} (2003), A new class of semiparametric estimators of the second order parameter,
{\em Port. Math.}, 60, 194--213.

\item{\sc Gadeikis, K. and V. Paulauskas} (2005), On the estimation of a change point in a tail index,
{\em Lith. Math. J.}, 45, 272--283.

\item{\sc Geluk, J.L. and L. Peng} (2000), An adaptive optimal estimate of the tail index for MA(1) time series,
{\em Statist. Probab. Lett.}, 46, 217--227.

\item{\sc Gomes, M.I., L. de Haan and  L. Henriques} (2008), Tail index estimation through the accommodation of bias in the weighted log-excesses,
{\em J. R. Stat. Soc. Ser. B. Stat. Methodol.}, 70, 31--53.

\item{\sc Gomes , M.I.  and M.J. Martins} (2002), Asymptotically unbiased estimators of the tail index
based on external estimation of the second order parameter,
{\em Extremes}, 5, 5--31.

\item{\sc Hall, P.} (1982), On some simple estimates of an exponent of regular variation,
{\em J. R. Stat. Soc. Ser. B. Stat. Methodol.}, 44, 37--42.

\item{\sc  Hill, B.M.} (1975), A simple general approach to inference about the tail of a distribution,
{\em Ann. Statist.}, 3, 1163--1174.

\item{\sc  LePage, R., M. Woodroofe and J. Zinn}
 (1981). Convergence to a stable distribution via order statistics. {\em Ann. Prob.} 9, 624--632.

\item{\sc  Li, J., Z. Peng and  S. Nadarajah} (2008), A class of unbiased location invariant Hill-type estimators for heavy tailed distributions,
{\em Electron. J. Stat.}, 2, 829--847.

\item{\sc Markovich, N.} (2007), {\em Nonparametric Analysis of Univariate Heavy-Tailed Data},
Jon Wiley \& Sons, Chichester.

\item{\sc  Nematollahi, A. R. and L. Tafakori} (2007),
On Comparison of the Tail Index of Heavy Tailed Distributions Using
Pitman's Measure of Closeness, {\em Appl. Math.\ Sci.},
1, 909--914.

\item{\sc  Paulauskas, V.} (2003), A New Estimator for a Tail Index,
{\em Acta Appl. Math. }, 79, 55--67.

\item{\sc  Paulauskas, V. and M. Vai\v ciulis } (2010), Some new
modifications of DPR estimator of the tail index, {\em preprint}.

\item{\sc Pickands, J.} (1975), Statistical inference using extreme order statistics,
{\em Ann. Statist. }, 3, 119--131.

\item{\sc Qi, Y.} (2010), On the tail index of a heavy tailed
distribution,{ \em Ann. Inst. Statist. Math.}, 62(2), 277--298.

\item{\sc Quintos, C., Zh. Fan and P. Phillips} (2001),Structural change tests in tail behavior and  the Asian crisis,
{\em Rev. Econom. Stud.}, 13, 633--663.

\item{\sc Resnick, S.  and C. Starica} (1997), Smoothing the Hill estimator,
{\em Adv. in Appl. Probab.}, 29, 271--293.

\item{\sc  Smith, R.L.} (1987), Estimating tails of probability distributions,
{\em Ann. Statist.}, 15, 1174--1207.

\item{\sc  Weissman, I.} (1978), Estimation of parameters and large quantiles based on the k largest observations,
{\em Journal of American Statistical Association }, 73, 812--815.
\end{description}

%\newpage
%\bibliographystyle{plain}
%\bibliography{probability7}

\end{document}